\documentclass[a4paper,11pt]{article}
\usepackage{amsmath}
\usepackage{amssymb}
\usepackage{tabularx}
\usepackage{enumerate}
\usepackage{setspace}
\usepackage{graphicx}
\usepackage{color}
\usepackage{graphicx}
\usepackage{caption}
\usepackage{subfigure}
\usepackage{float}

\newtheorem{theo}{Theorem}[section]
\newtheorem{prop}{Proposition}[section]
\newtheorem{lem}{Lemma}[section]
\newtheorem{rmk}{Remark}[section]

\newenvironment {Proof} {\noindent {\bf Proof.}}{\quad $\blacksquare$\par\vspace{3mm}}

\newcommand{\be}{\begin{equation}}
\newcommand{\ee}{\end{equation}}
\newcommand\bes{\begin{eqnarray}}
\newcommand\ees{\end{eqnarray}}
\newcommand{\bess}{\begin{eqnarray*}}
\newcommand{\eess}{\end{eqnarray*}}

\begin{document}
\setlength{\baselineskip}{15.2pt} \pagestyle{myheadings}

\title{ \bf  \LARGE Decay of solutions to anisotropic conservation laws with large initial data }
\date{\empty}
\author{ \large Kaiqiang Li,\ \ Weike Wang    \\
{ \normalsize School of Mathematical Sciences, Shanghai Jiao Tong University}\\
{ \normalsize Shanghai 200240, PR China}\\
{ \normalsize Email: kaiqiangli19@163.com(K.Li),\ \ wkwang@sjtu.edu.cn(W.Wang)} } \maketitle

\begin{quote}
\noindent {\bf Abstract:} {In this paper, we study the large time behavior of solutions to the Cauchy problem for the anisotropic conservation laws in two dimensional space. Without any smallness assumption on the initial data, the decay rates of solutions in $L^2$ space and homogeneous Sobolev space $\dot{H}^\gamma$ are obtained by using the method of time-frequency decomposition and the classical energy method. }


\end{quote}

\newcommand\HI{{\bf I}}
\section{Introduction.}
\setcounter{equation}{0}
In this paper, we study the large time behavior of solutions $u=u(x,t)$ to the Cauchy problem for the multi-dimensional anisotropic conservation laws with different fractional dissipation in diverse directions:
\begin{align}
 & u_t+\sum\limits_{i=1}^n\Lambda^{\alpha_i}_{x_i}u+\sum\limits_{i=1}^nf_i(u)_{x_i}=0, & (x,t)\in\mathbb{R}^n\times\mathbb{R}_+,\label{eq:1.1a} \\
 & u(x,0)=u_0(x), & x\in\mathbb{R}^n,\label{eq:1.2a}
\end{align}
where $\alpha_i\in(1,2]$, $f_i$'s are sufficiently smooth and satisfy
$$f_i(u)=O(|u|^{1+\kappa}),\ 1\leq i\leq n,\ 1\leq\kappa\in\mathbb{Z}^+.$$
In order to make sense of the effect about the diverse power in different directions clearly, we only consider the large time behavior of solutions for the above problem with $f(u)=u^2$ in two dimensional space:
\begin{align}
 & u_t+\Lambda^{\alpha_1}_xu+\Lambda_y^{\alpha_2}u+uu_x+uu_y=0, & (x,y,t)\in\mathbb{R}^2\times\mathbb{R}_+,\label{eq:1.1} \\
 & u(x,y,0)=u_0(x,y), & (x,y)\in\mathbb{R}^2,\label{eq:1.2}
\end{align}
where $\alpha_1,\alpha_2\in(1,2]$ and $\Lambda_x^{\alpha_1}, \Lambda_y^{\alpha_2}$ are the pseudo-differential operators defined via the Fourier transform
\begin{equation}\label{eq:1.3}
  \widehat{\Lambda_x^{\alpha_1}v}(\xi_1,\xi_2)=(\xi_1^2)^{\alpha_1/2}\widehat{v}(\xi_1,\xi_2)=|\xi_1|^{\alpha_1}\widehat{v}(\xi_1,\xi_2),
\end{equation}
and
\begin{equation}\label{eq:1.4}
  \widehat{\Lambda_y^{\alpha_2}v}(\xi_1,\xi_2)=(\xi_2^2)^{\alpha_2/2}\widehat{v}(\xi_1,\xi_2)=|\xi_2|^{\alpha_2}\widehat{v}(\xi_1,\xi_2).
\end{equation}

As the simplified prototype of more complicated problems arising in continuum mechanics, equation \eqref{eq:1.1} with $\alpha_1=\alpha_2=2$ has been widely considered for many years since it was studied by Burgers in the 1940s. [For $n=1$, the equality \eqref{eq:1.1} is just the well-known viscous Burgers equation.]

In one dimensional space case, it is well-known that the pure Burgers' equation can give rise to shocks even if the initial data is smooth enough. On one hand, when $\alpha_1=2$, it provides an accessible model for studying the interaction between nonlinear and dissipative phenomenon. In this case, Hopf \cite{EH} first studied the existence of solutions with nonlinear term $uu_x$ and then the optimal decay rates in $L^p$ norms and $H^m$ norms have been obtained by Kotlow \cite{DK} and Schonbek \cite{MS,MS1}, respectively. On the other hand, the value $\alpha_1=1$ is a threshold for the occurrence of singularity in finite time or the global regularity (see \cite{NJJ,HDD,JTJ,AFR}). Recently, Kiselev, Nazarov and Shterenberg \cite{AFR} considered in the circle $\mathbb{S}^1$ to show the finite time blow up for the supercritical case $0<\alpha<1$ and the global well-posedness in $H^{\frac{1}{2}}(\mathbb{S}^1)$ for the critical case $\alpha=1$ and in $H^s(\mathbb{S}^1)$($s>3/2-\alpha$) for the subcritical case $1<\alpha<2$. Also, Dong, Du and Li \cite{HDD} considered both of spaces $\mathbb{S}^1$ and $\mathbb{R}$ to show the finite time blow up for the supercritical case and the global well-posedness in $H^{\frac{1}{2}}(\mathbb{S}^1)$, $H^{\frac{1}{2}}(\mathbb{R})$ for the critical case and in $L^{1/(\alpha-1)}(\mathbb{S}^1)$, $L^{1/(\alpha-1)}(\mathbb{R})$ $(1<\alpha<2)$ for the subcritical case. Finite time blow up for the supercritical case is also shown by Alibaud, Droniou and Vovelle \cite{NJJ}. Miao and Wu \cite{MW} showed the global well-posedness in the critical case $\alpha=1$ for the initial data in the Besov spaces $\dot{B}_{p,1}^{\frac{1}{p}}(\mathbb{R})$ with $1\leq p<\infty$. On the large time behavior, Karch, Miao and Xu \cite{GCX} considered the subcritical case $1<\alpha<2$ in one dimensional space to show that the large time asymptotic is described by the rarefaction waves. More results can be found in \cite{NCG,BPGW,PGW,PTW1,AO,LAK,AK,AKN}.

In two dimensional space case,  Xin \cite{ZPX} has first investigated the stability of the planar rarefaction wave. Ito \cite{KI} has shown the convergence rate toward the planar rarefaction wave. However, in both papers, the smallness of initial disturbance is essentially assumed. Li \cite{FF} studied the time-asymptotic behavior of solutions to \eqref{eq:1.1} when $\alpha_1=\alpha_2\in(1,2]$, and obtained the optimal decay rates in $L^2$ and $\dot{H}^\gamma$ space for arbitrarily large initial data. Recently, Karch, Pudelko and Xu \cite{GAX} considered the two-dimensional convection-diffusion equation with $\alpha_1=\alpha_2=\alpha$ and showed that the large time behavior of solutions was described either by rarefaction waves, or diffusion waves, or suitable self-similar solutions, depending on the value of $\alpha$. We refer the reader to, e.g., \cite{DK,MK,MS,MS1} for an overview of known results and additional references.

To the author's best knowledge, this paper is the first trial to consider anisotropic conservation laws for \eqref{eq:1.1}-\eqref{eq:1.2}. In fact, the system can be taken as the classical fractional conservation laws in subcritical case due to the value of $\alpha_1, \alpha_2\in(1,2]$, so the global existence of solutions can be obtained by the standard method as the results for Burgers' equation in subcritical case, which the proof is given in Proposition \ref{prop:3.5} briefly. However, our aim is to make a detailed estimate for the large time behavior of the solutions with large initial data, and finally, we get the decay rates of the solutions in $L^2$ space and homogeneous Sobolev space $\dot{H}^{\gamma}$. Our result shows that the decay rates of solutions for anisotropic equation are determined by the different powers of diffusion terms in diverse directions. We first introduce the time-frequency cut-off operator which decompose the solutions into two parts: the low frequency part and the high frequency part. On the basis of the maximum principle, we get the decay rate of solutions in $L^2$ space by using of the energy method. Then, a preliminary decay rate in $\dot{H}^\gamma$ space is obtained by the energy method, which combining with Gagliardo-Nirenberg's inequality will give us a decay estimate in $L^\infty$ space. In the end, we can derive the decay rates of solutions in $\dot{H}^\gamma$ space on the basis of the $L^\infty$ estimate.

The paper is organized as follows. In Section 2, we introduce some notations and state our main result. In Section 3, we give some inequalities and draw into the maximum principle for equation \eqref{eq:1.1}. The decay rate of the solutions in $L^2$ space are obtained in Section 4. In Section 5, we first give some lemmas which are the tools for obtaining the main result, then we get the decay rate in $\dot{H}^\gamma$ space for the low frequency part and the high frequency part by different methods, respectively. In Section 6, we give two lemmas which could be used in the proof of the Theorem \ref{th:6.1}.

\section{Notations and main results.}
\setcounter{equation}{0}

We now list some notations that will be used in subsequent sections. As usual, the Fourier transform $\widehat{f}$ of $f$ is given by
$$\widehat{f}(\xi,t)\equiv(\mathcal{F}f)(\xi,t)=\int_{\mathbb{R}^n}f(x,t)e^{-ix\cdot\xi}dx.$$
The inverse Fourier transform is:
$$f(x,t)=(2\pi)^{-n}\int_{\mathbb{R}^n}\widehat{f}(\xi,t)e^{ix\cdot\xi}d\xi,\ \ \ x\in\mathbb{R}^n.$$
Denote $\Lambda=(-\Delta)^{\frac{1}{2}}$. Obviously
$$\widehat{\Lambda f}(\xi)=|\xi|\widehat{f}(\xi).$$
The space $\dot{W}^{\gamma,p}(\mathbb{R}^2)$ denotes the Sobolev space normed in two dimensional space by
$$\|\cdot\|_{\dot{W}^{\gamma,p}}=\|\Lambda^\gamma\cdot\|_{L^p},$$
where $\gamma\in\mathbb{R}$ and $\Lambda^\gamma$ is defined by Fourier transform
$$\widehat{\Lambda^\gamma f}(\xi)=|\xi|^{\gamma}\widehat{f}(\xi).$$
When $p=2$, we abbreviate $\dot{W}^{\gamma,2}$ to $\dot{H}^\gamma$, denoting the homogeneous Sobolev space.

Also, $\mathcal{D}$ is equivalent with $C_c^\infty$ which is the Fr$\acute{\mathrm{e}}$chet space of $C^\infty$ functions.

Throughout the paper, we use $\|\cdot\|_{L^p}$ to denote norm of $L^p(\mathbb{R}^2)$, and $C$ stands constant which its value may be different from line to line.

Now, let us describe our main theorem as follows.
\begin{theo}\label{th:2.1}
Let $\alpha_1, \alpha_2\in(1,2]$. If $u=u(x,y,t)$ is a solution for equation \eqref{eq:1.1} with initial data $u_0\in L^1(\mathbb{R}^2)\cap L^\infty(\mathbb{R}^2)$, then for any integer $\gamma\geq1$, we have

(1)\ the decay rate in $L^2$ space:
\begin{equation*}
  \|u(t,\cdot,\cdot)\|_{L^2}\leq C(1+t)^{-\frac{1}{2}\left(\frac{1}{\alpha_1}+\frac{1}{\alpha_2}\right)}.
\end{equation*}

(2)\ the decay rate in $\dot{H}^\gamma$ space:
\begin{equation*}
  \|u(t,\cdot,\cdot)\|_{\dot{H}^\gamma}\leq C(1+t)^{-\frac{1}{2}\left(\frac{1}{\alpha_1}+\frac{1}{\alpha_2}\right)-\frac{1}{2}\left(\frac{2\gamma+\lambda}{\alpha}-1\right)},
\end{equation*}
where $\alpha=\max\{\alpha_1,\alpha_2\}$, $\lambda=\min\{\alpha_1,\alpha_2\}$ and the constant $C$ only depends on the initial data $u_0$ and $\gamma$.
\end{theo}
\begin{rmk}\label{rm:2.1}
Especially, when $\alpha_1=\alpha_2=2$, we get the decay rates of solutions in $L^2$ space which is optimal and in homogeneous Sobolev space $\dot{H}^\gamma$, respectively:
\begin{equation*}
  \|u(t,\cdot,\cdot)\|_{L^2}\leq C(1+t)^{-\frac{1}{2}},
\end{equation*}
and
\begin{equation*}
  \|u(t,\cdot,\cdot)\|_{\dot{H}^\gamma}\leq C(1+t)^{-\frac{1+\gamma}{2}}.
\end{equation*}
Moreover, the above two estimates have been obtained in \cite{FF}.
\end{rmk}
\begin{rmk}\label{rm:2.2}
The result also holds for multi-dimensional space which the proof is more technique and complex in comparison with the case for two dimensional space. For convenience, we give the main result and its brief proof in Section 6.
\end{rmk}

\section{Preliminaries.}
\setcounter{equation}{0}
In this section, we recall some important inequalities and results for the proof of our main theorem.

\subsection{Some inequalities.}
To get the decay rate of the solutions for \eqref{eq:1.1} in homogeneous Sobolev space, we need to draw into the following inequalities.
\begin{lem}\label{le:3.1}(Gronwall's inequality \cite{CC})
Suppose that $a<b$ and let $\rho, \varphi$ and $\phi$ be non-negative continuous functions defined on the interval $[a,b]$. Moreover, suppose that $\rho$ is differentiable on $(a,b)$ with non-negative continuous derivative $\dot{\rho}$. If, for all $t\in[a,b]$,
$$\phi(t)\leq\rho(t)+\int_a^t\varphi(s)\phi(s)ds,$$
then
$$\phi(t)\leq\rho(t)e^{\int_a^t\varphi(s)ds}$$
for all $t\in[a,b]$.
\end{lem}
\begin{lem}\label{le:3.2}(Gagliardo-Nirenberg's inequality \cite{TC})
Let $1\leq p,q,r\leq\infty$ and let $j, m$ be two integers, $0\leq j<m$. If
$$\frac{1}{p}=\frac{j}{n}+a\left(\frac{1}{r}-\frac{m}{n}\right)+\frac{1-a}{q}$$
for some $a\in[\frac{j}{m},1]$ $(a<1\ if\ r>1\ and\ m-j-\frac{n}{r}=0)$, then there exists constant $C=C(n,m,j,a,q,r)$ such that
$$\sum\limits_{|\beta|=j}\|\nabla^\beta u\|_{L^p}\leq C\left(\sum\limits_{|\beta|=m}\|\nabla^\beta u\|_{L^r}\right)^a\|u\|_{L^q}^{1-a}$$
for every $u\in \mathcal{D}(\mathbb{R}^n)$.
\end{lem}
\begin{lem}\label{le:3.22}(Moser-type inequality)
Let $m\in\mathbb{N}$, then there is a constant $C=C(m,n)$ such that for all $f,g\in H^{m}(\mathbb{R}^n)\cap L^{\infty}(\mathbb{R}^n)$ and $\beta\in\mathbb{N}_0^n$, $|\beta|\leq m$, the following inequality holds:
$$\|\nabla^\beta(fg)\|_{L^2}\leq C\left(\|f\|_{L^\infty}\|\nabla^mg\|_{L^2}+\|g\|_{L^\infty}\|\nabla^mf\|_{L^2}\right).$$
\end{lem}
\subsection{Maximum principle.}
In this subsection, we describe the maximum principle in $L^p$ space, $1\leq p\leq\infty$, i.e., the $L^p$ norm of solutions for equation \eqref{eq:1.1} is controlled by the $L^p$ norm of the initial data.

The key point to prove the maximum principle is the positivity lemma as stated below. The original version was first presented in Resnick \cite{SR}. In \cite{PDJ}, Constantin, C$\acute{\mathrm{o}}$rdoba and Wu gave a detailed proof of the maximum principle for the Q-G equation. Our proof of the maximum principle is similar to that given in \cite{AD}, which is valid for more general flows.
\begin{lem}\label{le:3.3}(Positivity Lemma).
Let $0\leq\alpha\leq2$, $(x,y)\in\mathbb{R}^2$ and $u, \Lambda^{\alpha}u\in L^p(\mathbb{R}^2)$ with $1\leq p<\infty$, we have
$$\int_{\mathbb{R}^2}\Lambda^{\alpha}u|u|^{p-2}udxdy\geq0.$$
Especially, we can obtain
$$\int_{\mathbb{R}_y}\left[\int_{\mathbb{R}_x}\Lambda_x^{\alpha}u|u|^{p-2}udx\right]dy\geq0,\ \ \ \int_{\mathbb{R}_x}\left[\int_{\mathbb{R}_y}\Lambda_y^{\alpha}u|u|^{p-2}udy\right]dx\geq0.$$
\end{lem}
A detailed proof of the above lemma can be found in \cite{AD} for the Q-G equation, which also works for equation \eqref{eq:1.1}.
\begin{prop}\label{prop:3.4}(Maximum Principle).
Let $u$ be a smooth solution of equation \eqref{eq:1.1} with $\alpha_1, \alpha_2\in(1,2]$. Then for $1\leq p\leq\infty$, we have
\begin{equation}\label{eq:3.1}
  \|u\|_{L^p}\leq\|u_0\|_{L^p}.
\end{equation}
\end{prop}
\begin{Proof}
For $1\leq p<\infty$, using the positivity Lemma \ref{le:3.3}, we have
\begin{align}\label{eq:3.2}
  &\frac{d}{dt}\int_{\mathbb{R}^2}|u|^pdxdy\nonumber\\
= & p\int_{\mathbb{R}^2}|u|^{p-2}uu_tdxdy\nonumber\\
= & p\int_{\mathbb{R}^2}|u|^{p-2}u\left[-uu_x-uu_y-\Lambda_x^{\alpha_1}u-\Lambda_y^{\alpha_2}u\right]dxdy\nonumber\\
= & -p\int_{\mathbb{R}^2}|u|^{p-2}u\left(\Lambda_x^{\alpha_1}u+\Lambda_y^{\alpha_2}u\right)dxdy\nonumber\\
= & -p\int_{\mathbb{R}_y}\left[\int_{\mathbb{R}_x}|u|^{p-2}u\Lambda_x^{\alpha_1}udx\right]dy
    -p\int_{\mathbb{R}_x}\left[\int_{\mathbb{R}_y}|u|^{p-2}u\Lambda_y^{\alpha_2}udy\right]dx\nonumber\\
\leq & 0.
\end{align}
Thus, for $p\in[1,\infty)$, we obtain
$$\|u\|_{L^p}\leq\|u_0\|_{L^p}.$$
For $p=\infty$, the maximum principle is valid by \cite{JTJ}.
\end{Proof}

\subsection{Existence.}
In this subsection, we describe the global existence of solutions for the problem \eqref{eq:1.1}-\eqref{eq:1.2} by the classical method.
\begin{prop}\label{prop:3.5}(Existence).
Let $\alpha_1,\alpha_2\in(1,2]$, assume that the initial data $u_0\in L^1(\mathbb{R}^2)\cap H^\gamma(\mathbb{R}^2)$, $\gamma\geq2$, then there exists a global solution $u=u(t,x,y)$ for the problem \eqref{eq:1.1}-\eqref{eq:1.2} and it satisfies
\begin{equation*}
  u\in C([0,+\infty);H^\gamma(\mathbb{R}^2)),\ \gamma\geq2.
\end{equation*}
\end{prop}

Since $\alpha_1,\alpha_2\in(1,2],$ the equation \eqref{eq:1.1} can be taken as the classical fractal Burgers' equation in subcritical case. Then the global existence of solutions can be obtained basic on the results of the subcritical Burgers' equation or quasi-geostrophic equation. It could be proved by the standard method. Thus, we only give a framework for proof. The first step is to construct a local-in-time solution $u=u(t,x,y)$ on an interval $[0,T]$ for certain small $T>0$ which depends on the initial data $u_0$. Here, we should use the Banach contraction mapping principle. In the next step, by the maximum principle and some important inequalities we can improve the regularity of solutions. Then we can extend the local solution to be a global one. More details can be found in \cite{AD,PV,WW}.

\section{Decay rate in $L^2$.}
\setcounter{equation}{0}
In order to obtain the decay rate of solutions for equation \eqref{eq:1.1}-\eqref{eq:1.2} in $L^2$ norm, we first introduce the time-frequency cut-off operator which divides the solutions into two parts. Let
\begin{equation}\label{eq:4.1}
  \chi_0(\eta)=\left\{\aligned
  & 1,\ \ \ &\ \ \mathrm{if}\ \ |\eta|\leq1,\\
  & 0,\ \ \ &\ \ \mathrm{if}\ \ |\eta|>2,
  \endaligned\right.
\end{equation}
be a smooth function. Define the time-frequency cut-off operator $\chi(t,D)$ with the symbol $\chi(t,\xi)=\chi_0(\mu^{-1}(1+t)(|\xi_1|^{\alpha_1}+|\xi_2|^{\alpha_2}))$, where $\mu$ is a constant.

For a solution $u$ of \eqref{eq:1.1}, we can decompose it into two parts: the low frequency part $u_L$ and the high frequency part $u_H$, where $u_L(x,y,t)=\chi(t,D)u(x,y,t)$ and $u_H(x,y,t)=(1-\chi(t,D))u(x,y,t)$. Then $u_L$ and $u_H$ satisfy the following equations, respectively:
\begin{align}
  & \partial_t u_L+\Lambda_x^{\alpha_1}u_L+\Lambda_y^{\alpha_2}u_L+(uu_x)_L+(uu_y)_L=\partial_t\chi(t,D)u,\label{eq:4.2} \\
  & \partial_t u_H+\Lambda_x^{\alpha_1}u_H+\Lambda_y^{\alpha_2}u_H+(uu_x)_H+(uu_y)_H=-\partial_t\chi(t,D)u.\label{eq:4.3}
\end{align}

\subsection{Decay rate of low frequency part in $L^2$.}
Set $|\tau_1|^{\alpha_1}=\mu^{-1}(1+t)|\xi_1|^{\alpha_1}$, $|\tau_2|^{\alpha_2}=\mu^{-1}(1+t)|\xi_2|^{\alpha_2}$. According to the $L^2$-norm for the Fourier transform of $u_L$:
\begin{align}\label{eq:4.4}
          \|u_L(t,\cdot,\cdot)\|_{L^2}^2 = & \|\widehat{u_L}(t,\cdot,\cdot)\|_{L^2}^2
    =    \int_{\mathbb{R}^2}|\widehat{u}|^2(t,\xi_1,\xi_2)\chi_0^2\left(\mu^{-1}(1+t)(|\xi_1|^{\alpha_1}+|\xi_2|^{\alpha_2})\right)d\xi_1d\xi_2\nonumber  \\
\leq &   \|\widehat{u}(t,\cdot,\cdot)\|_{L^\infty}^2\int_{\mathbb{R}^2}\chi_0^2\left(\mu^{-1}(1+t)(|\xi_1|^{\alpha_1}+|\xi_2|^{\alpha_2})\right)
         d\xi_1d\xi_2\nonumber\\
\leq &   C\|u(t,\cdot,\cdot)\|_{L^1}^2\int_{\mathbb{R}^2}(1+t)^{-(\frac{1}{\alpha_1}+\frac{1}{\alpha_2})}
         \chi_0^2\left(|\tau_1|^{\alpha_1}+|\tau_2|^{\alpha_2}\right)d\tau_1d\tau_2\nonumber\\
\leq &   C(1+t)^{-(\frac{1}{\alpha_1}+\frac{1}{\alpha_2})}.
\end{align}
Thus we get the decay rate of the low frequency part $u_L$ in $L^2$:
\begin{equation}\label{eq:4.5}
  \|u_L(t,\cdot,\cdot)\|_{L^2}\leq C(1+t)^{-\frac{1}{2}(\frac{1}{\alpha_1}+\frac{1}{\alpha_2})}.
\end{equation}

\subsection{The decay rate of $u$ in $L^2$ space.}
We have the following interesting law between the decay rate of low frequency part $u_L$ and the decay rate of solution $u$.
\begin{prop}\label{prop:4.1}
If $\|u_L(t,\cdot,\cdot)\|_{L^2}\leq C(1+t)^{-\sigma}$, then
\begin{equation}\label{eq:4.6}
  \|u(t,\cdot,\cdot)\|_{L^2}\leq C(1+t)^{-\sigma}.
\end{equation}
\end{prop}
\begin{Proof}
Taking $L^2$ energy estimate of \eqref{eq:1.1}, we get
\begin{equation}\label{eq:4.7}
  \frac{1}{2}\frac{d}{dt}\|u\|_{L^2}^2+\|\Lambda_x^{\alpha_1/2}u\|_{L^2}^2+\|\Lambda_y^{\alpha_2/2}u\|_{L^2}^2=0,
\end{equation}
since $\int_{\mathbb{R}^2}\left(uu_x+uu_y\right)udxdy=0$ by the divergence theorem.

Let $\epsilon(t)=\mu(1+t)^{-1}$. We have
\begin{align}\label{eq:4.8}
 \|\Lambda_x^{\alpha_1/2}u\|_{L^2}^2+\|\Lambda_y^{\alpha_2/2}u\|_{L^2}^2 =  &  \int_{\mathbb{R}^2}\left(|\xi_1|^{\alpha_1}+|\xi_2|^{\alpha_2}\right)|\widehat{u}|^2d\xi_1d\xi_2\nonumber\\
 \geq  &  \int_{|\xi_1|^{\alpha_1}+|\xi_2|^{\alpha_2}\geq\epsilon(t)}\left(|\xi_1|^{\alpha_1}+|\xi_2|^{\alpha_2}\right)|\widehat{u}|^2d\xi_1d\xi_2\nonumber\\
 \geq  &  \epsilon(t)\int_{|\xi_1|^{\alpha_1}+|\xi_2|^{\alpha_2}\geq\epsilon(t)}|\widehat{u}|^2d\xi_1d\xi_2\nonumber\\
  =    &  \epsilon(t)\left[\|u\|_{L^2}^2-\int_{|\xi_1|^{\alpha_1}+|\xi_2|^{\alpha_2}\leq\epsilon(t)}|\widehat{u}|^2d\xi_1d\xi_2\right].
\end{align}
From \eqref{eq:4.7} and \eqref{eq:4.8}, we deduce that
\begin{align}\label{eq:4.9}
  & \frac{d}{dt}\|u\|_{L^2}^2+\epsilon(t)\|u\|_{L^2}^2 \nonumber\\
  \leq & \epsilon(t)\int_{|\xi_1|^{\alpha_1}+|\xi_2|^{\alpha_2}\leq\epsilon(t)}|\widehat{u}|^2d\xi_1d\xi_2\nonumber \\
\leq  & \epsilon(t)\int_{|\xi_1|^{\alpha_1}+|\xi_2|^{\alpha_2}\leq\epsilon(t)}
\chi_0^2\left(\mu^{-1}(1+t)\left(|\xi_1|^{\alpha_1}+|\xi_2|^{\alpha_2}\right)\right)|\widehat{u}|^2d\xi_1d\xi_2\nonumber\\
\leq & C\epsilon(t)\|u_L\|_{L^2}^2.
\end{align}
By the assumption, \eqref{eq:4.9} becomes
\begin{equation}\label{eq:4.10}
  \frac{d}{dt}\|u\|_{L^2}^2+\epsilon(t)\|u\|_{L^2}^2\leq C(1+t)^{-2\sigma-1}.
\end{equation}
Multiplying \eqref{eq:4.10} by $e^{\int_0^t\mu(1+\tau)^{-1}d\tau}=(1+t)^\mu$ and integrating from $0$ to $t$, we get
\begin{equation}\label{eq:4.11}
  (1+t)^\mu\|u\|_{L^2}^2\leq \|u_0\|_{L^2}^2+C(1+t)^{\mu-2\sigma}.
\end{equation}
Thus,
$$\|u\|_{L^2}^2\leq (1+t)^{-\mu}\|u_0\|_{L^2}^2+C(1+t)^{-2\sigma}.$$
Taking $\mu$ to be any constant greater that $2\sigma$, we obtain
$$\|u(t,\cdot,\cdot)\|_{L^2}\leq C(1+t)^{-\sigma}.$$
\end{Proof}
Combining \eqref{eq:4.5} and Proposition \ref{prop:4.1}, we can get the decay rate of $u$ in $L^2$ space:
\begin{equation}\label{eq:4.12}
  \|u(t,\cdot,\cdot)\|_{L^2}\leq C(1+t)^{-\frac{1}{2}\left(\frac{1}{\alpha_1}+\frac{1}{\alpha_2}\right)}.
\end{equation}

\section{Decay rate in $\dot{H}^{\gamma}$.}
\setcounter{equation}{0}
In this section, we get the decay rate of the solutions $u$ in homogeneous Sobolev space. We also decompose the solution into two parts: the low frequency part and the high frequency part which will use different methods. For the decay rate of the low frequency part we use the maximum principle, while for the high frequency part we use the energy estimate. Before giving the proof of our main result, we just draw into some lemmas which will be used in the following.

\subsection{Some lemmas.}
The following lemma is needed for the proof of Lemma \ref{le:5.2}. A similar estimate can be found in \cite{PG}.
\begin{lem}\label{le:5.1}
If $u_0\in L^1(\mathbb{R}^2)$, then the solution $u$ of \eqref{eq:1.1} satisfies the inequality
\begin{equation}\label{eq:5.1}
|\widehat{u}(t,\xi_1,\xi_2)|^2\leq\|u_0\|_{L^1}^2+C|\xi|\int_0^t\|u(\tau)\|_{L^2}^2\|u(\tau)\|_{L^1}d\tau.
\end{equation}
\end{lem}
\begin{Proof}
By Fourier transform, we can rewrite equation \eqref{eq:1.1} as
$$\partial_t\widehat{u}+\left(|\xi_1|^{\alpha_1}+|\xi_2|^{\alpha_2}\right)\widehat{u}+\left(\widehat{uu_x}+\widehat{uu_y}\right)=0.$$
Using
\begin{align*}
  \left|\widehat{uu_x}+\widehat{uu_y}\right|= & C\left|\int_{\mathbb{R}^2}e^{-(ix\xi_1+iy\xi_2)}\left(\partial_xu^2+\partial_yu^2\right)dxdy\right| \\
 \leq  & C|\xi|\int_{\mathbb{R}^2}|u^2|dxdy=C|\xi|\|u\|_{L^2}^2.
\end{align*}
We obtain that
\begin{align*}
\left|\widehat{u}(t,\xi_1,\xi_2)\right|^2\leq&\left|\widehat{u_0}(\xi_1,\xi_2)\right|^2
+C|\xi|\int_0^t\|u(\tau)\|_{L^2}^2|\widehat{u}(\tau,\xi_1,\xi_2)|d\tau\\
\leq & \|u_0\|_{L^1}^2+C|\xi|\int_0^t\|u(\tau)\|_{L^2}^2\|u(\tau)\|_{L^1}d\tau.
\end{align*}
So this complete the proof of this lemma.
\end{Proof}

\begin{lem}\label{le:5.2}
We have the following estimate:
\begin{equation}\label{eq:5.2}
  \int_0^t(1+\tau)^\mu\left(\|\Lambda_x^{\alpha_1/2}u\|_{L^2}^2+\|\Lambda_y^{\alpha_2/2}u\|_{L^2}^2\right)d\tau\leq C(1+t)^{\mu-\left(\frac{1}{\alpha_1}+\frac{1}{\alpha_2}\right)}
\end{equation}
where $C>0$ is a constant and $\mu>\frac{1}{\alpha_1}+\frac{1}{\alpha_2}$.
\end{lem}
\begin{Proof}
Assume $\epsilon(t)=\mu(1+t)^{-1}$. By the energy estimate \eqref{eq:4.7} we have
$$\frac{d}{dt}\|u\|_{L^2}^2+2\|\Lambda_x^{\alpha_1/2}u\|_{L^2}^2+2\|\Lambda_y^{\alpha_2/2}u\|_{L^2}^2=0.$$
Using the estimate \eqref{eq:4.8}, we obtain
\begin{align*}
\frac{d}{dt}\|u\|_{L^2}^2+\|\Lambda_x^{\alpha_1/2}u\|_{L^2}^2 & +\|\Lambda_y^{\alpha_2/2}u\|_{L^2}^2
   =   -\left(\|\Lambda_x^{\alpha_1/2}u\|_{L^2}^2+\|\Lambda_y^{\alpha_2/2}u\|_{L^2}^2\right)  \\
\leq  &  -\epsilon(t)\left[\|u\|_{L^2}^2-\int_{|\xi_1|^{\alpha_1}+|\xi_2|^{\alpha_2}\leq\epsilon(t)}|\widehat{u}|^2d\xi_1d\xi_2\right].
\end{align*}
Therefore, by Lemma \ref{le:5.1}, we get
\begin{align*}
     &  \frac{d}{dt}\|u\|_{L^2}^2+\|\Lambda_x^{\alpha_1/2}u\|_{L^2}^2+\|\Lambda_y^{\alpha_2/2}u\|_{L^2}^2+\epsilon(t)\|u\|_{L^2}^2  \\
\leq &  \epsilon(t)\int_{|\xi_1|^{\alpha_1}+|\xi_2|^{\alpha_2}\leq\epsilon(t)}|\widehat{u}|^2d\xi_1d\xi_2\\
\leq &  \epsilon(t)\int_{|\xi_1|^{\alpha_1}+|\xi_2|^{\alpha_2}\leq\epsilon(t)}
\left(\|u_0\|_{L^1}^2+C|\xi|\int_0^t\|u(\tau)\|_{L^2}^2\|u(\tau)\|_{L^1}d\tau\right)d\xi_1d\xi_2\\
\leq & C(1+t)^{-1-\left(\frac{1}{\alpha_1}+\frac{1}{\alpha_2}\right)}.
\end{align*}
As in the proof of Proposition \ref{prop:4.1}, we have
\begin{align*}
  &(1+t)^\mu\|u\|_{L^2}^2+\int_0^t(1+\tau)^\mu\left(\|\Lambda_x^{\alpha_1/2}u\|_{L^2}^2+\|\Lambda_y^{\alpha_2/2}u\|_{L^2}^2\right)d\tau \\
  \leq & \|u_0\|_{L^2}^2+C\int_0^t(1+\tau)^{\mu-1-\left(\frac{1}{\alpha_1}+\frac{1}{\alpha_2}\right)}d\tau \\
  \leq & C(1+t)^{\mu-\left(\frac{1}{\alpha_1}+\frac{1}{\alpha_2}\right)}.
\end{align*}
Taking $\mu$ to be any constant greater than $\frac{1}{\alpha_1}+\frac{1}{\alpha_2}$, then the inequality \eqref{eq:5.2} now follow from \eqref{eq:4.12}.
\end{Proof}
\begin{lem}\label{le:5.3}
Let $\alpha_1,\alpha_2\in(1,2]$, $\gamma\geq1$ is any integer. Then there exists constant $C>0$ such that
\begin{equation}\label{eq:5.3}
  \|\nabla^\gamma\Lambda_x^{1-\frac{\alpha_1}{2}}u\|_{L^2}\leq C\left(\|\nabla^\gamma\Lambda_x^{{\alpha_1}/{2}}u\|_{L^2}^{\theta_1}\|\Lambda_x^{{\alpha_1}/{2}}u\|_{L^2}^{1-\theta_1}
  +\|\nabla^\gamma\Lambda_y^{{\alpha_2}/{2}}u\|_{L^2}^{\theta_2}\|\Lambda_y^{{\alpha_2}/{2}}u\|_{L^2}^{1-\theta_2}\right)
\end{equation}
where
\begin{equation}\label{eq:5.4}
  \theta_1=\frac{\gamma+1-\alpha_1}{\gamma},\ \ \ \ \ \ \theta_2=\frac{2\gamma+2-\alpha_1-\alpha_2}{2\gamma}.
\end{equation}
Similarly, we also have
\begin{equation}\label{eq:5.5}
  \|\nabla^\gamma\Lambda_y^{1-\frac{\alpha_2}{2}}u\|_{L^2}\leq C\left(\|\nabla^\gamma\Lambda_x^{{\alpha_1}/{2}}u\|_{L^2}^{\theta_3}\|\Lambda_x^{{\alpha_1}/{2}}u\|_{L^2}^{1-\theta_3}
  +\|\nabla^\gamma\Lambda_y^{{\alpha_2}/{2}}u\|_{L^2}^{\theta_4}\|\Lambda_y^{{\alpha_2}/{2}}u\|_{L^2}^{1-\theta_4}\right)
\end{equation}
where
\begin{equation}\label{eq:5.6}
\theta_3=\frac{2\gamma+2-\alpha_1-\alpha_2}{2\gamma},\ \ \ \ \ \ \theta_4=\frac{\gamma+1-\alpha_2}{\gamma}.
\end{equation}
\end{lem}
\begin{Proof}
In order to prove the inequality \eqref{eq:5.3}, we only need to prove the following two inequalities:
\begin{equation}\label{eq:5.7}
  \|\nabla_x^{\gamma+1-\frac{\alpha_1}{2}}u\|_{L^2}\leq C\|\nabla^\gamma\Lambda_x^{\alpha_1/2}u\|_{L^2}^{\theta_1}\|\Lambda_x^{\alpha_1/2}u\|_{L^2}^{1-\theta_1},
\end{equation}
and
\begin{equation}\label{eq:5.8}
  \|\nabla_y^{\gamma+1-\frac{\alpha_1}{2}}u\|_{L^2}\leq C\|\nabla^\gamma\Lambda_y^{\alpha_2/2}u\|_{L^2}^{\theta_2}\|\Lambda_y^{\alpha_2/2}u\|_{L^2}^{1-\theta_2},
\end{equation}
since
\begin{equation}\label{eq:5.9}
  \|\nabla^\gamma\Lambda_x^{1-\frac{\alpha_1}{2}}u\|_{L^2}\leq\|\nabla^{\gamma+1-\frac{\alpha_1}{2}}u\|_{L^2}\leq C\left(\|\nabla_x^{\gamma+1-\frac{\alpha_1}{2}}u\|_{L^2}+ \|\nabla_y^{\gamma+1-\frac{\alpha_1}{2}}u\|_{L^2}\right),
\end{equation}
and $\theta_1, \theta_2$ are described in \eqref{eq:5.4}.

Let us recall this short and elegant reasoning. For every $R>0$, we decompose the $L^2$-norm of the Fourier transform as follows:
\begin{align}\label{eq:5.10}
  \|\nabla_x^{\gamma+1-\frac{\alpha_1}{2}}u\|_{L^2}^2 = & C\int_{\mathbb{R}^2}|\xi_1|^{2\gamma+2-\alpha_1}|\widehat{u}|^2d\xi_1d\xi_2 \nonumber\\
 =  & C\int_{|\xi|\leq R}|\xi_1|^{\alpha_1}|\widehat{u}|^2|\xi_1|^{2\gamma+2-2\alpha_1}d\xi+C\int_{|\xi|\geq R}|\xi_1|^{2\gamma+2-\alpha_1}|\widehat{u}|^2d\xi\nonumber\\
 \leq & CR^{2\gamma+2-2\alpha_1}\int_{\mathbb{R}^2}|\xi_1|^{\alpha_1}|\widehat{u}|^2d\xi+C\int_{|\xi|\geq R}|\xi_1|^{2\gamma+2-\alpha_1}|\widehat{u}|^2d\xi\nonumber\\
 \leq & CR^{2\gamma+2-2\alpha_1}\|\Lambda_x^{\alpha_1/2}u\|_{L^2}^2+C\int_{|\xi|\geq R}|\xi_1|^{2\gamma+2-\alpha_1}|\widehat{u}|^2d\xi,
\end{align}
where $|\xi|^2=|\xi_1|^2+|\xi_2|^2$.

For the second term in the right hand side of \eqref{eq:5.10}, it can be viewed as
\begin{align}\label{eq:5.11}
   & \int_{|\xi|\geq R}|\xi_1|^{2\gamma+2-\alpha_1}|\widehat{u}|^2d\xi\nonumber \\
=  & \left(\int_{|\xi_1|\geq|\xi_2|,|\xi|\geq R}+\int_{|\xi_2|\geq|\xi_1|,|\xi|\geq R}\right)\left(|\xi_1|^{2\gamma+2-\alpha_1}|\widehat{u}|^2\right)d\xi\nonumber\\
\leq &C\int_{|\xi_1|\geq\frac{1}{\sqrt{2}}R}|\xi_1|^{2\gamma}|\xi_1|^{\alpha_1}|\xi_1|^{2-2\alpha_1}|\widehat{u}|^2d\xi
+C\int_{|\xi_2|\geq\frac{1}{\sqrt{2}}R}|\xi_1|^{\alpha_1}|\xi_1|^{2\gamma+2-2\alpha_1}|\widehat{u}|^2d\xi\nonumber\\
\leq& CR^{2-2\alpha_1}\|\nabla^\gamma\Lambda_x^{\alpha_1/2}u\|_{L^2}^2
+C\int_{|\xi_2|\geq\frac{1}{\sqrt{2}}R}|\xi_1|^{\alpha_1}|\xi_2|^{2\gamma}|\xi_2|^{2-2\alpha_1}|\widehat{u}|^2d\xi\nonumber\\
\leq & CR^{2-2\alpha_1}\|\nabla^\gamma\Lambda_x^{\alpha_1/2}u\|_{L^2}^2,
\end{align}
where we use the fact that
\begin{equation}\label{eq:5.12}
  |\xi_1|^{2}=C\left(|\xi_1|^2+|\xi_2|^2\right),\ \ \ \ \ when\ \ |\xi_1|\geq|\xi_2|,
\end{equation}
and
\begin{equation}\label{eq:5.13}
  |\xi_2|^{2}=C\left(|\xi_1|^2+|\xi_2|^2\right),\ \ \ \ \ when\ \ |\xi_2|\geq|\xi_1|.
\end{equation}
Combining the inequality \eqref{eq:5.10} with \eqref{eq:5.11}, we have
\begin{equation}\label{eq:5.14}
   \|\nabla_x^{\gamma+1-\frac{\alpha_1}{2}}u\|_{L^2}^2\leq CR^{2\gamma+2-2\alpha_1}\|\Lambda_x^{\alpha_1/2}u\|_{L^2}^2+CR^{2-2\alpha_1}\|\nabla^\gamma\Lambda_x^{\alpha_1/2}u\|_{L^2}^2.
\end{equation}
It suffices to choose
$$R=\left(\frac{\|\nabla^\gamma\Lambda_x^{\alpha_1/2}u\|_{L^2}}{\|\Lambda_x^{\alpha_1/2}u\|_{L^2}}\right)^{\frac{1}{\gamma}}$$
for obtaining inequality \eqref{eq:5.7}, i.e.,
$$\|\nabla_x^{\gamma+1-\frac{\alpha_1}{2}}u\|_{L^2}
\leq C\|\nabla^\gamma\Lambda_x^{\alpha_1/2}u\|_{L^2}^{\theta_1}\|\Lambda_x^{\alpha_1/2}u\|_{L^2}^{1-\theta_1},$$
where $\theta_1=\frac{\gamma+1-\alpha_1}{\gamma}$.

Similarly,
\begin{align}\label{eq:5.15}
  \|\nabla_y^{\gamma+1-\frac{\alpha_1}{2}}u\|_{L^2}^2 = & C\int_{\mathbb{R}^2}|\xi_2|^{2\gamma+2-\alpha_1}|\widehat{u}|^2d\xi_1d\xi_2 \nonumber\\
 =  & C\int_{|\xi|\leq R}|\xi_2|^{\alpha_2}|\widehat{u}|^2|\xi_2|^{2\gamma+2-\alpha_1-\alpha_2}d\xi+C\int_{|\xi|\geq R}|\xi_2|^{2\gamma+2-\alpha_1}|\widehat{u}|^2d\xi\nonumber\\
 \leq & CR^{2\gamma+2-\alpha_1-\alpha_2}\int_{\mathbb{R}^2}|\xi_2|^{\alpha_2}|\widehat{u}|^2d\xi+C\int_{|\xi|\geq R}|\xi_2|^{2\gamma+2-\alpha_1}|\widehat{u}|^2d\xi\nonumber\\
 \leq & CR^{2\gamma+2-\alpha_1-\alpha_2}\|\Lambda_y^{\alpha_2/2}u\|_{L^2}^2+C\int_{|\xi|\geq R}|\xi_2|^{2\gamma+2-\alpha_1}|\widehat{u}|^2d\xi,
\end{align}
For the second term in the right hand side of \eqref{eq:5.15}, using equality \eqref{eq:5.12} and \eqref{eq:5.13}, we have
\begin{align}\label{eq:5.16}
   & \int_{|\xi|\geq R}|\xi_2|^{2\gamma+2-\alpha_1}|\widehat{u}|^2d\xi\nonumber \\
=  & \left(\int_{|\xi_1|\geq|\xi_2|,|\xi|\geq R}+\int_{|\xi_2|\geq|\xi_1|,|\xi|\geq R}\right)\left(|\xi_2|^{2\gamma+2-\alpha_1}|\widehat{u}|^2\right)d\xi\nonumber\\
\leq & C\int_{|\xi_1|\geq\frac{1}{\sqrt{2}}R}|\xi_2|^{\alpha_2}|\xi_1|^{2\gamma+2-\alpha_1-\alpha_2}|\widehat{u}|^2d\xi
+C\int_{|\xi_2|\geq\frac{1}{\sqrt{2}}R}|\xi_2|^{\alpha_2}|\xi_2|^{2\gamma+2-\alpha_1-\alpha_2}|\widehat{u}|^2d\xi\nonumber\\
\leq& CR^{2-\alpha_1-\alpha_2}\|\nabla^\gamma\Lambda_y^{\alpha_2/2}u\|_{L^2}^2
+C\int_{|\xi_2|\geq\frac{1}{\sqrt{2}}R}|\xi_2|^{\alpha_2}|\xi_2|^{2\gamma}|\xi_2|^{2-\alpha_1-\alpha_2}|\widehat{u}|^2d\xi\nonumber\\
\leq & CR^{2-\alpha_1-\alpha_2}\|\nabla^\gamma\Lambda_y^{\alpha_2/2}u\|_{L^2}^2.
\end{align}
Combining inequality \eqref{eq:5.15} with \eqref{eq:5.16}, we obtain that
\begin{equation}\label{eq:5.17}
   \|\nabla_y^{\gamma+1-\frac{\alpha_1}{2}}u\|_{L^2}^2\leq CR^{2\gamma+2-\alpha_1-\alpha_2}\|\Lambda_y^{\alpha_2/2}u\|_{L^2}^2+CR^{2-\alpha_1-\alpha_2}\|\nabla^\gamma\Lambda_y^{\alpha_2/2}u\|_{L^2}^2.
\end{equation}
It suffices to choose
$$R=\left(\frac{\|\nabla^\gamma\Lambda_y^{\alpha_2/2}u\|_{L^2}}{\|\Lambda_y^{\alpha_2/2}u\|_{L^2}}\right)^{\frac{1}{\gamma}}$$
to obtain the inequality \eqref{eq:5.8}, i.e.,
$$\|\nabla_y^{\gamma+1-\frac{\alpha_1}{2}}u\|_{L^2}
\leq C\|\nabla^\gamma\Lambda_y^{\alpha_2/2}u\|_{L^2}^{\theta_2}\|\Lambda_y^{\alpha_2/2}u\|_{L^2}^{1-\theta_2},$$
where $\theta_2=\frac{2\gamma+2-\alpha_1-\alpha_2}{2\gamma}$.

Similarly, the inequality \eqref{eq:5.5} can also be obtained by the same way which described as above. Thus, Lemma \ref{le:5.3} is finished.
\end{Proof}

\begin{lem}\label{le:5.4}
Let $\alpha_1,\alpha_2\in(1,2]$, and $\gamma\geq1$ is any integer. Then there is constant $C>0$ satisfy
\begin{equation}\label{eq:5.18}
  \|\nabla^\gamma\Lambda_x^{1-\frac{\alpha_1}{2}}u\|_{L^2}\leq C\left(\|\nabla^\gamma\Lambda_x^{\alpha_1/2}u\|_{L^2}^{s_1}\|\nabla^\gamma u\|_{L^2}^{1-s_1}
  +\|\nabla^\gamma\Lambda_y^{\alpha_2/2}u\|_{L^2}^{s_2}\|\nabla^\gamma u\|_{L^2}^{1-s_2}\right)
\end{equation}
where
\begin{equation}\label{eq:5.19}
  s_1=\frac{2-\alpha_1}{\alpha_1},\ \ \ \ \ \ \ s_2=\frac{2-\alpha_1}{\alpha_2}.
\end{equation}
Similarly, we also have
\begin{equation}\label{eq:5.20}
  \|\nabla^\gamma\Lambda_y^{1-\frac{\alpha_2}{2}}u\|_{L^2}\leq C\left(\|\nabla^\gamma\Lambda_x^{\alpha_1/2}u\|_{L^2}^{s_3}\|\nabla^\gamma u\|_{L^2}^{1-s_3}
  +\|\nabla^\gamma\Lambda_y^{\alpha_2/2}u\|_{L^2}^{s_4}\|\nabla^\gamma u\|_{L^2}^{1-s_4}\right)
\end{equation}
where
\begin{equation}\label{eq:5.21}
s_3=\frac{2-\alpha_2}{\alpha_1},\ \ \ \ \ \ \ s_4=\frac{2-\alpha_2}{\alpha_2}.
\end{equation}
\end{lem}
\begin{Proof}
Here, we only prove the first inequality \eqref{eq:5.18} since inequality \eqref{eq:5.20} can be obtained with the similar method. As we can see that
\begin{equation}\label{eq:5.22}
  \|\nabla^\gamma\Lambda_x^{1-\frac{\alpha_1}{2}}u\|_{L^2}\leq \|\nabla^{\gamma+1-\frac{\alpha_1}{2}}u\|_{L^2}\leq C\left(\|\nabla_x^{\gamma+1-\frac{\alpha_1}{2}}u\|_{L^2}+\|\nabla_y^{\gamma+1-\frac{\alpha_1}{2}}u\|_{L^2}\right).
\end{equation}
So, we only need to prove the following two inequalities:
\begin{equation}\label{eq:5.23}
  \|\nabla_x^{\gamma+1-\frac{\alpha_1}{2}}u\|_{L^2}\leq C\|\nabla^\gamma\Lambda_x^{\alpha_1/2}u\|_{L^2}^{s_1}\|\nabla^\gamma u\|_{L^2}^{1-s_1},
\end{equation}
and
\begin{equation}\label{eq:5.24}
  \|\nabla_y^{\gamma+1-\frac{\alpha_1}{2}}u\|_{L^2}\leq C\|\nabla^\gamma\Lambda_y^{\alpha_2/2}u\|_{L^2}^{s_2}\|\nabla^\gamma u\|_{L^2}^{1-s_2}.
\end{equation}
For every $R>0$, we decompose the $L^2$-norm of the Fourier transform as follows:
\begin{align}\label{eq:5.25}
  \|\nabla_x^{\gamma+1-\frac{\alpha_1}{2}}u\|_{L^2}^2 = & C\int_{\mathbb{R}^2}|\xi_1|^{2\gamma+2-\alpha_1}|\widehat{u}|^2d\xi_1d\xi_2 \nonumber\\
 =  & C\int_{|\xi|\leq R}|\xi_1|^{2\gamma}|\widehat{u}|^2|\xi_1|^{2-\alpha_1}d\xi+C\int_{|\xi|\geq R}|\xi_1|^{2\gamma+2-\alpha_1}|\widehat{u}|^2d\xi\nonumber\\
 \leq & CR^{2-\alpha_1}\int_{\mathbb{R}^2}|\xi|^{2\gamma}|\widehat{u}|^2d\xi+C\int_{|\xi|\geq R}|\xi_1|^{2\gamma+2-\alpha_1}|\widehat{u}|^2d\xi\nonumber\\
 \leq & CR^{2-\alpha_1}\|\nabla^\gamma u\|_{L^2}^2+C\int_{|\xi|\geq R}|\xi_1|^{2\gamma+2-\alpha_1}|\widehat{u}|^2d\xi.
\end{align}
For the second term on the right hand side of \eqref{eq:5.25}, we have
\begin{align}\label{eq:5.26}
   & \int_{|\xi|\geq R}|\xi_1|^{2\gamma+2-\alpha_1}|\widehat{u}|^2d\xi\nonumber \\
=  & \left(\int_{|\xi_1|\geq|\xi_2|,|\xi|\geq R}+\int_{|\xi_2|\geq|\xi_1|,|\xi|\geq R}\right)\left(|\xi_1|^{2\gamma+2-\alpha_1}|\widehat{u}|^2\right)d\xi\nonumber\\
\leq &C\int_{|\xi_1|\geq\frac{1}{\sqrt{2}}R}|\xi_1|^{2\gamma}|\xi_1|^{\alpha_1}|\xi_1|^{2-2\alpha_1}|\widehat{u}|^2d\xi
+C\int_{|\xi_2|\geq\frac{1}{\sqrt{2}}R}|\xi_1|^{\alpha_1}|\xi_1|^{2\gamma+2-2\alpha_1}|\widehat{u}|^2d\xi\nonumber\\
\leq& CR^{2-2\alpha_1}\|\nabla^\gamma\Lambda_x^{\alpha_1/2}u\|_{L^2}^2
+C\int_{|\xi_2|\geq\frac{1}{\sqrt{2}}R}|\xi_1|^{\alpha_1}|\xi_2|^{2\gamma}|\xi_2|^{2-2\alpha_1}|\widehat{u}|^2d\xi\nonumber\\
\leq & CR^{2-2\alpha_1}\|\nabla^\gamma\Lambda_x^{\alpha_1/2}u\|_{L^2}^2.
\end{align}
Combining the above two inequalities \eqref{eq:5.25} and \eqref{eq:5.26}, we obtain
\begin{equation*}
   \|\nabla_x^{\gamma+1-\frac{\alpha_1}{2}}u\|_{L^2}^2\leq CR^{2-\alpha_1}\|\nabla^\gamma u\|_{L^2}^2+CR^{2-2\alpha_1}\|\nabla^\gamma\Lambda_x^{\alpha_1/2}u\|_{L^2}^2.
\end{equation*}
It suffices to choose $R=\left(\|\nabla^\gamma\Lambda_x^{\alpha_1/2}u\|_{L^2}^2/\|\nabla^\gamma u\|_{L^2}^2\right)^{1/{\alpha_1}}$ to get the inequality \eqref{eq:5.23}, i.e.,
\begin{equation*}
  \|\nabla_x^{\gamma+1-\frac{\alpha_1}{2}}u\|_{L^2}
\leq C\|\nabla^\gamma\Lambda_x^{\alpha_1/2}u\|_{L^2}^{s_1}\|\nabla^\gamma u\|_{L^2}^{1-s_1},
\end{equation*}
where $s_1=\frac{2-\alpha_1}{\alpha_1}$.

Similarly, we have
\begin{equation*}
  \|\nabla_y^{\gamma+1-\frac{\alpha_1}{2}}u\|_{L^2}^2\leq CR^{2-\alpha_1}\|\nabla^\gamma u\|_{L^2}^2+CR^{2-\alpha_1-\alpha_2}\|\nabla^\gamma\Lambda_y^{\alpha_2/2}u\|_{L^2}^2.
\end{equation*}
Choosing $R=\left(\|\nabla^\gamma\Lambda_y^{\alpha_2/2}u\|_{L^2}^2/\|\nabla^\gamma u\|_{L^2}^2\right)^{1/{\alpha_2}}$, we can get the inequality \eqref{eq:5.24}, i.e.,
\begin{equation*}
  \|\nabla_y^{\gamma+1-\frac{\alpha_1}{2}}u\|_{L^2}
\leq C\|\nabla^\gamma\Lambda_y^{\alpha_2/2}u\|_{L^2}^{s_2}\|\nabla^\gamma u\|_{L^2}^{1-s_2},
\end{equation*}
where $s_2=\frac{2-\alpha_1}{\alpha_2}$.

Combining the inequality \eqref{eq:5.23} with \eqref{eq:5.24} then we finish the proof of the inequality \eqref{eq:5.18}.

At the same time, we can get the inequality \eqref{eq:5.20} by the same way. Thus, the proof of Lemma \ref{le:5.4} is finished.
\end{Proof}

\subsection{Decay rate of low frequency part in $\dot{H}^{\gamma}$.}
As before, we set $|\eta_1|^{\alpha_1}=\mu^{-1}(1+t)|\xi_1|^{\alpha_1}$, $|\eta_2|^{\alpha_2}=\mu^{-1}(1+t)|\xi_2|^{\alpha_2}$, we have
\begin{align*}
     & \|\nabla^\gamma u_L(t,\cdot,\cdot)\|_{L^2}^2 =  \||\xi|^{\gamma}\widehat{u_L}(t,\cdot,\cdot)\|_{L^2}^2 \\
  =  & \int_{\mathbb{R}^2}|\xi|^{2\gamma}|\widehat{u}|^2\chi_0^2\left(\mu^{-1}(1+t)(|\xi_1|^{\alpha_1}+|\xi_2|^{\alpha_2})\right)d\xi\\
\leq & \|\widehat{u}\|_{L^\infty}^2\int_{\mathbb{R}^2}|\xi|^{2\gamma}\chi_0^2\left(\mu^{-1}(1+t)(|\xi_1|^{\alpha_1}+|\xi_2|^{\alpha_2})\right)d\xi\\
\leq & \|u\|_{L^1}^2\int_{\mathbb{R}^2}\left(|\xi_1|^2+|\xi_2|^2\right)^\gamma\chi_0^2\left(\mu^{-1}(1+t)(|\xi_1|^{\alpha_1}+|\xi_2|^{\alpha_2})\right)d\xi\\
\leq & C(1+t)^{-\left(\frac{1}{\alpha_1}+\frac{1}{\alpha_2}\right)}\int_{\mathbb{R}^2}\left((1+t)^{-\frac{2\gamma}{\alpha_1}}|\eta_1|^{2\gamma}
+(1+t)^{-\frac{2\gamma}{\alpha_2}}|\eta_2|^{2\gamma}\right)\chi_0^2\left(|\eta_1|^{\alpha_1}+|\eta_2|^{\alpha_2}\right)d\eta_1d\eta_2\\
\leq & C(1+t)^{-\left(\frac{1}{\alpha_1}+\frac{1}{\alpha_2}\right)}\left[(1+t)^{-\frac{2\gamma}{\alpha_1}}+(1+t)^{-\frac{2\gamma}{\alpha_2}}\right]\\
\leq & C(1+t)^{-\left(\frac{1}{\alpha_1}+\frac{1}{\alpha_2}\right)-\frac{2\gamma}{\alpha}},
\end{align*}
where $\alpha=\max\{\alpha_1,\alpha_2\}$.

Thus, we get the decay rate of low frequency part in $\dot{H}^\gamma$:
\begin{equation}\label{eq:5.27}
  \|u_L(t,\cdot,\cdot)\|_{\dot{H}^\gamma}\leq C(1+t)^{-\frac{1}{2}\left(\frac{1}{\alpha_1}+\frac{1}{\alpha_2}\right)-\frac{\gamma}{\alpha}}
\end{equation}
where $\alpha=\max\{\alpha_1,\alpha_2\}$.
\subsection{Decay rate of high frequency part in $\dot{H}^{\gamma}$.}
In what follows, we will obtain the decay rate for the solutions in homogeneous Sobolev space $\dot{H}^\gamma$. First, we obtain a preliminary decay rate by use of the maximum principle and Lemma \ref{le:5.2}, then combining with Gagliardo-Nirenberg's inequality will give us a decay estimate of the solutions in $L^\infty$ space. Second, we can obtain the decay rate of solutions for \eqref{eq:1.1}-\eqref{eq:1.2} in $\dot{H}^\gamma$ space on the basis of the $L^\infty$ estimate.

We know that $u_H$ satisfies the equation \eqref{eq:4.3}:
\begin{equation*}
  \partial_t u_H+\Lambda_x^{\alpha_1}u_H+\Lambda_y^{\alpha_2}u_H+(uu_x)_H+(uu_y)_H=-\partial_t\chi(t,D)u.
\end{equation*}
In order to estimate $u_H$, multiplying \eqref{eq:4.3} by $\nabla^{2\gamma}u_H$ and integrating it in $\mathbb{R}^2$, we get
\begin{align}\label{eq:5.28}
    & \frac{1}{2}\frac{d}{dt}\|\nabla^\gamma u_H\|_{L^2}^2+\|\nabla^\gamma\Lambda_x^{\alpha_1/2}u_H\|_{L^2}^2+\|\nabla^\gamma\Lambda_y^{\alpha_2/2}u_H\|_{L^2}^2 \nonumber\\
 =  & -\int_{\mathbb{R}^2}(uu_x)_H\nabla^{2\gamma}u_Hdxdy-\int_{\mathbb{R}^2}(uu_y)_H\nabla^{2\gamma}u_Hdxdy-\int_{\mathbb{R}^2}\partial_t\chi(t,D)u\nabla^{2\gamma}u_Hdxdy\nonumber\\
 := & J_1+J_2+J_3.
\end{align}
For the nonlinear term $J_1$, by H$\ddot{\mathrm{o}}$lder inequality, Lemma \ref{le:3.3} and \ref{le:5.3}, we have
\begin{align}\label{eq:5.29}
  &|J_1|\leq C\left|\int_{\mathbb{R}^2}\nabla^\gamma\partial_x(u^2)\nabla^{\gamma}udxdy\right|\leq C\|u\|_{L^\infty}\left|\int_{\mathbb{R}^2}\nabla^\gamma\Lambda_x^{\alpha_1/2}u\nabla^\gamma\Lambda_x^{1-\alpha_1/2}udxdy\right|\nonumber\\
  \leq & C\|u\|_{L^\infty}\|\nabla^\gamma\Lambda_x^{\alpha_1/2}u\|_{L^2}\|\nabla^\gamma\Lambda_x^{1-\alpha_1/2}u\|_{L^2}\nonumber\\
  \leq & C\|\nabla^\gamma\Lambda_x^{\alpha_1/2}u\|_{L^2}\left(\|\nabla^\gamma\Lambda_x^{\alpha_1/2}u\|_{L^2}^{\theta_1}\|\Lambda_x^{\alpha_1/2}u\|_{L^2}^{1-\theta_1}
  +\|\nabla^\gamma\Lambda_y^{\alpha_2/2}u\|_{L^2}^{\theta_2}\|\Lambda_y^{\alpha_2/2}u\|_{L^2}^{1-\theta_2}\right)\nonumber\\
  \leq & C\|\nabla^\gamma\Lambda_x^{\alpha_1/2}u\|_{L^2}^{1+\theta_1}\|\Lambda_x^{\alpha_1/2}u\|_{L^2}^{1-\theta_1}
  + C\|\nabla^\gamma\Lambda_x^{\alpha_1/2}u\|_{L^2}\|\nabla^\gamma\Lambda_y^{\alpha_2/2}u\|_{L^2}^{\theta_2}\|\Lambda_y^{\alpha_2/2}u\|_{L^2}^{1-\theta_2}\nonumber\\
  \leq & \frac{1}{4}\|\nabla^\gamma\Lambda_x^{\alpha_1/2}u\|_{L^2}^2+\frac{1}{4}\|\nabla^\gamma\Lambda_y^{\alpha_2/2}u\|_{L^2}^2
  +C\left(\|\Lambda_x^{\alpha_1/2}u\|_{L^2}^2+\|\Lambda_y^{\alpha_2/2}u\|_{L^2}^2\right),
\end{align}
where
\begin{equation*}
  \theta_1=\frac{\gamma+1-\alpha_1}{\gamma},\ \ \ \ \ \ \theta_2=\frac{2\gamma+2-\alpha_1-\alpha_2}{2\gamma}.
\end{equation*}
Similarly, for the second term on the right hand side of \eqref{eq:5.28}, we have
\begin{align}\label{eq:5.30}
  &|J_2|\leq C\left|\int_{\mathbb{R}^2}\nabla^\gamma\partial_y(u^2)\nabla^{\gamma}udxdy\right|\leq C\|u\|_{L^\infty}\left|\int_{\mathbb{R}^2}\nabla^\gamma\Lambda_y^{\alpha_2/2}u\nabla^\gamma\Lambda_y^{1-\alpha_2/2}udxdy\right|\nonumber\\
  \leq & \frac{1}{4}\|\nabla^\gamma\Lambda_x^{\alpha_1/2}u\|_{L^2}^2+\frac{1}{4}\|\nabla^\gamma\Lambda_y^{\alpha_2/2}u\|_{L^2}^2
  +C\left(\|\Lambda_x^{\alpha_1/2}u\|_{L^2}^2+\|\Lambda_y^{\alpha_2/2}u\|_{L^2}^2\right).
\end{align}
For $J_3$, let $|\eta_1|^{\alpha_1}=\mu^{-1}(1+t)|\xi_1|^{\alpha_1}$, $|\eta_2|^{\alpha_2}=\mu^{-1}(1+t)|\xi_2|^{\alpha_2}$. We also assume $\rho=\mu^{-1}(1+t)\left(|\xi_1|^{\alpha_1}+|\xi_2|^{\alpha_2}\right)$ for convenience. Then
\begin{align}\label{eq:5.31}
  |J_3|
\leq &  C\int_{\mathbb{R}^2}|\widehat{u}|^2\mu^{-1}|\xi|^{2\gamma}\left(|\xi_1|^{\alpha_1}+|\xi_2|^{\alpha_2}\right)
\chi_0'(\rho)\left|1-\chi_0\left(\rho\right)\right|d\xi\nonumber\\
\leq &  C\|\widehat{u}\|_{L^\infty}^2\int_{\mathbb{R}^2}(1+t)^{-1}|\xi|^{2\gamma}\chi_0'(\rho)|1-\chi_0(\rho)|d\xi\nonumber\\
\leq & C\|u\|_{L^1}^2\int_{\mathbb{R}^2}(1+t)^{-1-\left(\frac{1}{\alpha_1}+\frac{1}{\alpha_2}\right)-\frac{2\gamma}{\alpha}}
\left(|\eta_1|^{2\gamma}+|\eta_2|^{2\gamma}\right)\chi_0'(\omega)\left|1-\chi_0(\omega)\right|d\eta_1d\eta_2\nonumber\\
\leq & C(1+t)^{-1-\left(\frac{1}{\alpha_1}+\frac{1}{\alpha_2}\right)-\frac{2\gamma}{\alpha}},
\end{align}
where $\alpha=\max\{\alpha_1,\alpha_2\}$ and we denote $\omega=|\eta_1|^{\alpha_1}+|\eta_2|^{\alpha_2}$.

Combining \eqref{eq:5.28}-\eqref{eq:5.31}, we obtain
\begin{align}\label{eq:5.32}
&  \frac{d}{dt}\|\nabla^\gamma u_H\|_{L^2}^2+\|\nabla^\gamma\Lambda_x^{\alpha_1/2}u_H\|_{L^2}^2+\|\nabla^\gamma\Lambda_y^{\alpha_2/2}u_H\|_{L^2}^2 \nonumber\\
\leq & \|\nabla^\gamma\Lambda_x^{\alpha_1/2}u_L\|_{L^2}^2+\|\nabla^\gamma\Lambda_y^{\alpha_2/2}u_L\|_{L^2}^2
+C\left(\|\Lambda_x^{\alpha_1/2}u\|_{L^2}^2+\|\Lambda_y^{\alpha_2/2}u\|_{L^2}^2\right)\nonumber\\
  &  +C(1+t)^{-1-\left(\frac{1}{\alpha_1}+\frac{1}{\alpha_2}\right)-\frac{2\gamma}{\alpha}}\nonumber\\
\leq & C(1+t)^{-\left(\frac{1}{\alpha_1}+\frac{1}{\alpha_2}\right)-\frac{2\gamma}{\alpha}}
\left((1+t)^{-1}+(1+t)^{-\frac{\lambda}{\alpha}}\right) +C\left(\|\Lambda_x^{\alpha_1/2}u\|_{L^2}^2+\|\Lambda_y^{\alpha_2/2}u\|_{L^2}^2\right)\nonumber\\
\leq & C(1+t)^{-\left(\frac{1}{\alpha_1}+\frac{1}{\alpha_2}\right)-\frac{2\gamma+\lambda}{\alpha}}
+C\left(\|\Lambda_x^{\alpha_1/2}u\|_{L^2}^2+\|\Lambda_y^{\alpha_2/2}u\|_{L^2}^2\right),
\end{align}
where $\alpha=\max\left\{\alpha_1,\alpha_2\right\}$ and $\lambda=\min\{\alpha_1,\alpha_2\}$.

In fact, using inequality \eqref{eq:5.27} we have
\begin{equation}\label{eq:5.33}
  \|\nabla^\gamma\Lambda_x^{\alpha_1/2}u_L\|_{L^2}^2\leq C\|\nabla^{\gamma+\frac{\alpha_1}{2}}u_L\|_{L^2}^2\leq C(1+t)^{-\left(\frac{1}{\alpha_1}+\frac{1}{\alpha_2}\right)-\frac{2\gamma+\alpha_1}{\alpha}},
\end{equation}
and
\begin{equation}\label{eq:5.34}
  \|\nabla^\gamma\Lambda_y^{\alpha_2/2}u_L\|_{L^2}^2\leq C\|\nabla^{\gamma+\frac{\alpha_2}{2}}u_L\|_{L^2}^2\leq C(1+t)^{-\left(\frac{1}{\alpha_1}+\frac{1}{\alpha_2}\right)-\frac{2\gamma+\alpha_2}{\alpha}}.
\end{equation}
According to the definition of the cut-off operator, we know that
\begin{align}\label{eq:5.35}
\|\nabla^\gamma\Lambda_x^{\alpha_1/2}u_H\|_{L^2}^2+\|\nabla^\gamma\Lambda_y^{\alpha_2/2}u_H\|_{L^2}^2
= & \int_{\mathbb{R}^2}|\xi|^{2\gamma}|\widehat{u}|^2\left(|\xi_1|^{\alpha_1}+|\xi_2|^{\alpha_2}\right)d\xi\nonumber\\
\geq & \mu(1+t)^{-1}\|\nabla^\gamma u_H\|_{L^2}^2.
\end{align}
Then we have
\begin{align}\label{eq:5.36}
\frac{d}{dt}\|\nabla^\gamma u_H\|_{L^2}^2+\mu(1+t)^{-1}\|\nabla^\gamma u_H\|_{L^2}^2\leq & C(1+t)^{-\left(\frac{1}{\alpha_1}+\frac{1}{\alpha_2}\right)-\frac{2\gamma+\lambda}{\alpha}}\nonumber\\
&+C\left(\|\Lambda_x^{\alpha_1/2}u\|_{L^2}^2+\|\Lambda_y^{\alpha_2/2}u\|_{L^2}^2\right).
\end{align}
Using Lemma \ref{le:5.2} and multiplying \eqref{eq:5.36} by $e^{\int_0^t\mu(1+\tau)^{-1}d\tau}=(1+t)^\mu$ and integrating from 0 to $t$, we have
\begin{align}\label{eq:5.37}
 (1+t)^\mu\|\nabla^\gamma u_H\|_{L^2}^2\leq & \|\nabla^\gamma u_{0H}\|_{L^2}^2
 +C\int_0^t(1+\tau)^\mu(1+\tau)^{-\left(\frac{1}{\alpha_1}+\frac{1}{\alpha_2}\right)-\frac{2\gamma+\lambda}{\alpha}}d\tau\nonumber\\
  & +C\int_0^t(1+\tau)^\mu\left(\|\Lambda_x^{\alpha_1/2}u\|_{L^2}^2+\|\Lambda_y^{\alpha_2/2}u\|_{L^2}^2\right)d\tau\nonumber\\
\leq & \|\nabla^\gamma u_{0H}\|_{L^2}^2+(1+t)^{\mu+1-\left(\frac{1}{\alpha_1}+\frac{1}{\alpha_2}\right)-\frac{2\gamma+\lambda}{\alpha}}\nonumber\\
&+C(1+t)^{\mu-\left(\frac{1}{\alpha_1}+\frac{1}{\alpha_2}\right)}\nonumber\\
\leq & \|\nabla^\gamma u_{0H}\|_{L^2}^2+C(1+t)^{\mu-\left(\frac{1}{\alpha_1}+\frac{1}{\alpha_2}\right)},
\end{align}
for the fact that $1-\frac{2\gamma+\lambda}{\alpha}<0$.

Taking $\mu$ to be any constant greater than $\frac{1}{\alpha_1}+\frac{1}{\alpha_2}$, we then get the decay rate of high frequency part in $\dot{H}^\gamma$:
\begin{equation}\label{eq:5.38}
  \|u_H(t,\cdot,\cdot)\|_{\dot{H}^{\gamma}}\leq C(1+t)^{-\frac{1}{2}\left(\frac{1}{\alpha_1}+\frac{1}{\alpha_2}\right)}.
\end{equation}
Combining \eqref{eq:5.27} with \eqref{eq:5.38}, we have
\begin{equation}\label{eq:5.39}
  \|\nabla^\gamma u\|_{L^2}\leq \|\nabla^\gamma u_H\|_{L^2}+\|\nabla^\gamma u_L\|_{L^2}\leq C(1+t)^{-\frac{1}{2}\left(\frac{1}{\alpha_1}+\frac{1}{\alpha_2}\right)}.
\end{equation}
This implies that when $\gamma=2$, we obtain
\begin{equation}\label{eq:5.40}
  \|u(t,\cdot,\cdot)\|_{\dot{H}^2}\leq C(1+t)^{-\frac{1}{2}\left(\frac{1}{\alpha_1}+\frac{1}{\alpha_2}\right)}.
\end{equation}
By Gagliardo-Nirenberg's inequality, we have
\begin{equation}\label{eq:5.41}
  \|u\|_{L^\infty}\leq C\|u\|_{\dot{H}^2}^{\frac{1}{2}}\|u\|_{L^2}^{\frac{1}{2}}\leq C(1+t)^{-\frac{1}{2}\left(\frac{1}{\alpha_1}+\frac{1}{\alpha_2}\right)}.
\end{equation}
Using Lemma \ref{le:5.4} and \eqref{eq:5.29}, \eqref{eq:5.30} and Young inequality, we can get
\begin{align}\label{eq:5.42}
 &      |J_1|+|J_2|\nonumber\\
\leq & C\|u\|_{L^\infty}\left(\|\nabla^\gamma\Lambda_x^{\alpha_1/2}u\|_{L^2}^{1+s_1}\|\nabla^\gamma u\|_{L^2}^{1-s_1}+\|\nabla^\gamma\Lambda_x^{\alpha_1/2}u\|_{L^2}\|\nabla^\gamma u\|_{L^2}^{1-s_2}\|\nabla^\gamma\Lambda_y^{\alpha_2/2}\|_{L^2}^{s_2}\right)\nonumber\\
  &  +C\|u\|_{L^\infty}\left(\|\nabla^\gamma\Lambda_y^{\alpha_2/2}u\|_{L^2}\|\nabla^\gamma\Lambda_x^{\alpha_1/2}u\|_{L^2}^{s_3}\|\nabla^\gamma u\|_{L^2}^{1-s_3}+\|\nabla^\gamma\Lambda_y^{\alpha_2/2}\|_{L^2}^{1+s_4}\|\nabla^\gamma u\|_{L^2}^{1-s_4}\right)\nonumber\\
\leq & C\left(\|u\|_{L^\infty}^{\frac{2}{1-s_1}}+\|u\|_{L^\infty}^{\frac{2}{1-s_2}}+\|u\|_{L^\infty}^{\frac{2}{1-s_3}}+\|u\|_{L^\infty}^{\frac{2}{1-s_4}}\right)
\|\nabla^\gamma u\|_{L^2}^2+\frac{1}{2}\|\nabla^\gamma\Lambda_x^{\alpha_1/2}u\|_{L^2}^2\nonumber\\
& +\frac{1}{2}\|\nabla^\gamma\Lambda_y^{\alpha_2/2}u\|_{L^2}^2\nonumber\\
\leq & \frac{1}{2}\|\nabla^\gamma\Lambda_x^{\alpha_1/2}u\|_{L^2}^2+\frac{1}{2}\|\nabla^\gamma\Lambda_y^{\alpha_2/2}u\|_{L^2}^2
+C\|u\|_{L^\infty}^2\|\nabla^\gamma u\|_{L^2}^2.
\end{align}
According to the inequality \eqref{eq:5.31}, we have
\begin{equation*}
  |J_3|\leq C(1+t)^{-1-\left(\frac{1}{\alpha_1}+\frac{1}{\alpha_2}\right)-\frac{2\gamma}{\alpha}}.
\end{equation*}
Similar with inequality \eqref{eq:5.32}, we then have
\begin{align}\label{eq:5.43}
      & \frac{d}{dt}\|\nabla^\gamma u_H\|_{L^2}^2+\|\nabla^\gamma\Lambda_x^{\alpha_1/2}u_H\|_{L^2}^2+\|\nabla^\gamma\Lambda_y^{\alpha_2/2} u_H\|_{L^2}^2\nonumber\\
 \leq & \|\nabla^\gamma\Lambda_x^{\alpha_1/2}u_L\|_{L^2}^2+\|\nabla^\gamma\Lambda_y^{\alpha_2/2}u_L\|_{L^2}^2+C\|u\|_{L^\infty}^2\|\nabla^\gamma u\|_{L^2}^2 +C(1+t)^{-1-\left(\frac{1}{\alpha_1}+\frac{1}{\alpha_2}\right)-\frac{2\gamma}{\alpha}}\nonumber\\
\leq & C(1+t)^{-\left(\frac{1}{\alpha_1}+\frac{1}{\alpha_2}\right)
-\frac{2\gamma}{\alpha}}\left((1+t)^{-1}+(1+t)^{-\frac{\lambda}{\alpha}}\right)
+C\|u\|_{L^\infty}^2\|\nabla^\gamma u\|_{L^2}^2\nonumber\\
\leq & C(1+t)^{-\left(\frac{1}{\alpha_1}+\frac{1}{\alpha_2}\right)-\frac{2\gamma+\lambda}{\alpha}}
+C(1+t)^{-\left(\frac{1}{\alpha_1}+\frac{1}{\alpha_2}\right)}\|\nabla^\gamma u\|_{L^2}^2,
\end{align}
where $\alpha=\max\left\{\alpha_1,\alpha_2\right\}$ and $\lambda=\min\{\alpha_1,\alpha_2\}$.

Multiplying the above inequality by $e^{\int_0^t\mu(1+\tau)^{-1}d\tau}=(1+t)^\mu$, combining \eqref{eq:5.27} with \eqref{eq:5.35} and then integrating from 0 to $t$, we get
\begin{align}
     & (1+t)^\mu\|\nabla^\gamma u_H\|_{L^2}^2\nonumber\\
\leq & \|\nabla^\gamma u_{0H}\|_{L^2}^2+C(1+t)^{\mu+1-\left(\frac{1}{\alpha_1}+\frac{1}{\alpha_2}\right)-\frac{2\gamma+\lambda}{\alpha}}\nonumber\\
     & +C\int_0^t(1+\tau)^{\mu-\left(\frac{1}{\alpha_1}+\frac{1}{\alpha_2}\right)}\left(\|\nabla^\gamma u_L\|_{L^2}^2+\|\nabla^\gamma u_H\|_{L^2}^2\right)d\tau\label{eq:5.44}\\
\leq & \|\nabla^\gamma u_{0H}\|_{L^2}^2+C(1+t)^{\mu+1-\left(\frac{1}{\alpha_1}+\frac{1}{\alpha_2}\right)-\frac{2\gamma+\lambda}{\alpha}}
+C\int_0^t(1+\tau)^{\mu-\left(\frac{1}{\alpha_1}+\frac{1}{\alpha_2}\right)}\|\nabla^\gamma u_H\|_{L^2}^2d\tau.\nonumber
\end{align}
By Gronwall's inequality, we have
\begin{align}\label{eq:5.45}
  (1+t)^\mu\|\nabla^\gamma u_H\|_{L^2}^2 \leq & C\left(1+(1+t)^{\mu+1-\left(\frac{1}{\alpha_1}+\frac{1}{\alpha_2}\right)-\frac{2\gamma+\lambda}{\alpha}}\right)
  e^{\int_0^t(1+\tau)^{-\left(\frac{1}{\alpha_1}+\frac{1}{\alpha_2}\right)}d\tau}\nonumber\\
  \leq & C+C(1+t)^{\mu+1-\left(\frac{1}{\alpha_1}+\frac{1}{\alpha_2}\right)-\frac{2\gamma+\lambda}{\alpha}}.
\end{align}
Taking $\mu$ to be any constant greater than $\frac{1}{\alpha_1}+\frac{1}{\alpha_2}+\frac{2\gamma+\lambda}{\alpha}-1$, we then obtain
\begin{equation}\label{eq:5.46}
  \|u_H(t,\cdot,\cdot)\|_{\dot{H}^\gamma}\leq C(1+t)^{\frac{1}{2}-\frac{1}{2}\left(\frac{1}{\alpha_1}+\frac{1}{\alpha_2}\right)-\frac{2\gamma+\lambda}{2\alpha}},
\end{equation}
where $\alpha=\max\{\alpha_1,\alpha_2\}$ and $\lambda=\min\{\alpha_1,\alpha_2\}$.

Finally, combining \eqref{eq:5.27} and \eqref{eq:5.46}, we conclude that there is a constant $C$ which only depends on $\gamma, \alpha_1, \alpha_2$ and the initial data $u_0$ such that
\begin{equation}\label{eq:5.47}
   \|u(t,\cdot,\cdot)\|_{\dot{H}^\gamma}\leq C(1+t)^{-\frac{1}{2}\left(\frac{1}{\alpha_1}+\frac{1}{\alpha_2}\right)-\frac{1}{2}\left(\frac{2\gamma+\lambda}{\alpha}-1\right)},
\end{equation}
where $\alpha=\max\{\alpha_1,\alpha_2\}$ and $\lambda=\min\{\alpha_1,\alpha_2\}$.

\section{The general case.}
\setcounter{equation}{0}
For the multi-dimensional space case, i.e., the Cauchy problem as follows:
\begin{align}
   & u_t+\sum\limits_{i=1}^n\Lambda_{x_i}^{\alpha_i}u+\sum\limits_{i=1}^nf_i(u)_{x_i}=0, & (x,t)\in\mathbb{R}^n\times\mathbb{R}_+,\label{eq:6.1} \\
   & u(x,0)=u_0(x),& x\in\mathbb{R}^n,\label{eq:6.2}
\end{align}
where the initial data $u_0$ is arbitrarily large and $\alpha_i\in(1,2]$, $\Lambda_{x_i}^{\alpha_i}$ is the pseudo-differential operator defined via the Fourier transform
\begin{equation}\label{eq:6.3}
  \widehat{\Lambda_{x_i}^{\alpha_i}v}(\xi)=|\xi_i|^{\alpha_i}\widehat{v}(\xi),
\end{equation}
with $|\xi|=\left(\sum\limits_{i=1}^n|\xi_i|^2\right)^{\frac{1}{2}}$. Here, $f_i$ is sufficiently smooth and satisfies
\begin{equation}\label{eq:6.4}
  f_i(u)=O(|u|^{1+\kappa}),\ \ 1\leq i\leq n,
\end{equation}
with $1\leq\kappa\in\mathbb{Z}^+$.

First, the global existence of solutions for the above problem can be obtained by using the method in \cite{AD,PV,WW}.

Then, on the basis of the global existence of solutions such as Proposition \ref{prop:3.5}, we give the main result which is similar with Theorem \ref{th:2.1}.
\begin{theo}\label{th:6.1}
Let $\alpha_i\in(1,2]$. If $u=u(x,t)$ is a solution of equation \eqref{eq:6.1} with the initial data $u_0\in L^1(\mathbb{R}^n)\cap L^\infty(\mathbb{R}^n)$, then for any integer $\gamma\geq1$, there exists a constant $C>0$, we have

(1)\ the decay rate in $L^2$ space:
\begin{equation*}
  \|u(t,\cdot)\|_{L^2}\leq C(1+t)^{-\frac{1}{2}\sum\limits_{i=1}^n\frac{1}{\alpha_i}}.
\end{equation*}

(2)\ the decay rate in $\dot{H}^\gamma$ space:
\begin{equation*}
  \|u(t,\cdot)\|_{\dot{H}^\gamma}\leq C(1+t)^{-\frac{1}{2}\sum\limits_{i=1}^n\frac{1}{\alpha_i}-\frac{1}{2}\left(\frac{2\gamma+\lambda}{\alpha}-1\right)}
\end{equation*}
where $\alpha=\max\{\alpha_1,\alpha_2,\cdot,\cdot,\cdot,\alpha_n\}$ and $\lambda=\min\{\alpha_1,\alpha_2,\cdot,\cdot,\cdot,\alpha_n\}$.
\end{theo}
Now, we draw into some lemmas which play an important role in our proof of the decay rate about the solutions for \eqref{eq:6.1}-\eqref{eq:6.2}.
\begin{lem}\label{le:6.1}
Assume that $f(u)$ is smooth enough, and $f(u)=O(|u|^{1+\kappa})$ when $|u|\leq\mu_0$, where $\kappa\geq1$ is an integer. For each integer $\gamma\geq0$, if $u$ satisfies $\|u\|_{L^\infty}\leq\mu_0$, then
\begin{equation}\label{eq:6.5}
  \|f(u)\|_{\dot{H}^\gamma}\leq C\|u\|_{\dot{H}^\gamma}\|u\|_{L^\infty}^\kappa,
\end{equation}
where $C$ is a constant depending on $\gamma$ and $\mu_0$.
\end{lem}
\begin{lem}\label{le:6.2}
Let $\alpha_i\in(1,2]$, $\gamma\geq1$ is an integer. Then there exists a constant $C>0$ such that
\begin{equation}\label{eq:6.6}
  \|\nabla^\gamma\Lambda_{x_i}^{1-\frac{\alpha_i}{2}}u\|_{L^2}\leq C\sum\limits_{i=1}^n\|\nabla^\gamma\Lambda_{x_i}^{\alpha_i/2}u\|_{L^2}^{\theta_i}\|\Lambda_{x_i}^{\alpha_i/2}u\|_{L^2}^{1-\theta_i},
\end{equation}
and
\begin{equation}\label{eq:6.7}
  \|\nabla^\gamma\Lambda_{x_i}^{1-\frac{\alpha_i}{2}}u\|_{L^2}\leq C\sum\limits_{i=1}^n\|\nabla^\gamma\Lambda_{x_i}^{\alpha_i/2}u\|_{L^2}^{s_i}\|\nabla^\gamma u\|_{L^2}^{1-s_i}.
\end{equation}
where $\theta_i$, $s_i$ depend on $\gamma$ and $\alpha_i$ and $\theta_i, s_i\in[0,1)$.
\end{lem}
The above lemma can be got such as the proof of Lemma \ref{le:5.3} and \ref{le:5.4} which we omit it for convenience.

According to Lemma \ref{le:6.1} and \ref{le:6.2}, using the same method which described above we can easily get the decay rate of solutions for problem \eqref{eq:6.1}-\eqref{eq:6.2}. In this section, we omit the detail proof of Theorem \ref{th:6.1} for convenience.
\section*{Acknowledgments}

The paper is supported by National Nature Science Foundation of China (No.11771284 and No.11871335). The first author would like to thank Prof. Xiongfeng Yang for his useful discussions.

\end{document}